\documentclass[10pt,twoside]{article}
\usepackage{graphicx}
\usepackage{amsmath}
\usepackage{Latex-document}

\title{\bf Face Numbers \vskip -2mm
of 4-Polytopes and 3-Spheres\vskip 6mm}

\author{G\"unter M. Ziegler\vspace*{-0.5cm}\thanks{Institute of
    Mathematics, MA 6-2, Technical University of Berlin, D-10623 Berlin, \newline
    Germany. E-mail: ziegler@math.tu-berlin.de,
    http://www.math.tu-berlin.de/\char126ziegler}}
\date{\vspace{-8mm}}
\usepackage{url,amssymb,paralist,color}
\newcommand\Q{{\mathbb Q}}
\newcommand\R{{\mathbb R}}
\newcommand\Z{{\mathbb Z}}
\newcommand\fatness{\mathrm{F}}
\newcommand\complexity{\mathrm{C}}
\newcommand\conv{\mathop\mathrm{conv}}
\newcommand\flag{\mathop\mathrm{flag}}
\newcommand\fourflag{\mathop\mathrm{flag}}
\newcommand\fvectors{\mbox{\rm{\it f}\/-Vectors}}
\newcommand\flagvectors{\mbox{\rm Flag-Vectors}}
\newcommand\closedfcone{\mbox{\rm$f$-Cone}}
\newcommand\closedflagcone{\mbox{\rm Flag-Cone}}
\newcommand\PP{{\cal P}_4}
\newcommand\EL{{\cal L}_4}
\newcommand\SP{{\cal S}_4}
\newcommand\PR{{\cal P}_4^{\Q}}
\newcommand\para[1]{\smallskip\noindent{\bf#1}}
\begin{document}
\maketitle

\thispagestyle{first} \setcounter{page}{625}


\begin{abstract}\vskip 3mm
Steinitz (1906) gave a remarkably simple
and explicit description of the set of all $f$-vectors
$f(P)=(f_0,f_1,f_2)$ of all $3$-dimensional convex polytopes.
His result also identifies the simple and the
simplicial $3$-dimensional polytopes as the only extreme cases.
Moreover, it can be extended to
strongly regular CW $2$-spheres (topological objects), and
further to Eulerian lattices of length~$4$ (combinatorial objects).

The analogous problems ``one dimension higher,'' about the
$f$-vectors and flag-vectors of $4$-dimensional convex polytopes
and their generalizations, are by far not solved, yet.  However,
the known facts already show that the answers will be much more
complicated than for Steinitz' problem.  In this lecture, we will
summarize the current state of knowledge.  We will put forward two
crucial parameters of \emph{fatness} and \emph{complexity}:
Fatness $\fatness(P):=\frac{{\displaystyle
f_1+f_2-20}}{{\displaystyle f_0+f_3-10}}$ is large if there are
many more edges and $2$-faces than there are vertices and facets,
while complexity $\complexity(P):=\frac{{\displaystyle
f_{03}-20}}{{\displaystyle f_0+f_3-10}}$ is large if every facet
has many vertices, and every vertex is in many facets. Recent
results suggest that these parameters might allow one to
differentiate between the cones of $f$- or flag-vectors~of
\begin{compactitem}[~$\bullet$]
\item connected Eulerian lattices of length~$5$ (combinatorial objects),
\item strongly regular CW $3$-spheres (topological objects),
\item convex $4$-polytopes (discrete geometric objects), and
\item rational convex $4$-polytopes (whose study involves arithmetic aspects).
\end{compactitem}
Further progress will depend on the derivation of tighter $f$-vector
inequalities for convex $4$-polytopes.
On the other hand, we will need new construction
methods that produce interesting polytopes which are far
from being simplicial or simple --- for example, very ``fat'' or
``complex'' $4$-polytopes. In this direction, I will report
about constructions (from joint work with Michael Joswig,
David Eppstein and Greg Kuperberg) that yield
\begin{compactitem}[~$\bullet$]
\item  strongly regular CW $3$-spheres of arbitrarily large fatness,
\item  convex $4$-polytopes of fatness larger than $5.048$, and
\item  rational convex $4$-polytopes of fatness larger than $5-\varepsilon$.
\end{compactitem}
\vskip 4.5mm


\end{abstract}

\vskip 12mm


\section{Introduction} 
\setzero\vskip-5mm \hspace{5mm }

Our knowledge about the combinatorics and geometry
of $4$-dimen\-sional convex polytopes is quite incomplete.
This assessment may come as a surprise: After all,
\begin{compactitem}[~$\bullet$]
\item $3$-dimensional polytopes have been objects of intensive study
  since antiquity,
\item properties of convex polytopes are essential to the
  geometry of Euclidean spaces,
\item the \emph{regular} polytopes (in all dimensions) were
  classified by Schl\"afli in 1850-52~\cite{Schla},
  exactly 150 years ago, and
\item modern polytope theory has achieved truly impressive results, in
  particular since the publication of Gr\"unbaum's
  volume~\cite{Gr1} in 1967, thirty-five years ago.
\end{compactitem}
Moreover, we have a rather satisfactory picture of $3$-dimensional
convex polytopes by now, where the essential combinatorial and
geometric properties of $3$-dimensional polytopes were isolated a long
time ago.  Here we mention the three results that will be most
relevant to our subsequent discussion:
\begin{compactenum}[~1.]
\item
  Steinitz~\cite{Stei3} characterized
  the $f$-vectors $(f_0,f_1,f_2)\in\Z^3$ of the $3$-dimensional
  convex polytopes: They are the integer points in the $2$-dimensional
  convex polyhedral cone that is defined by the three conditions
\[
f_0-f_1+f_2\ =\ 2,\qquad
f_2-4\ \le\ 2(f_0-4),\qquad
f_0-4\ \le\ 2(f_2-4).
\]
In particular, the $f$-vectors of the $3$-dimensional polytopes are
given as \emph{all} the integer points in a rational polyhedral cone.
Furthermore, Steinitz' result includes the characterization of the
polytopes with extremal $f$-vectors: The first inequality is
tight if and only if the polytope is simplicial, while the second
one is tight if and only if the polytope is simple.
\item In 1922, Steinitz~\cite{Stei1} published a characterization of the
  graphs of $3$-polytopes: They are all the planar,
  three-connected graphs on at least $4$ nodes.\\
  In modern terms (as reviewed below) and after some additional
  arguments the Steinitz theorem may be phrased as follows: Every
  connected finite Eulerian lattice of length~$4$ is the face lattice
  of a rational convex $3$-polytope.
\item The famous Koebe--Andreev--Thurston circle-packing
  theorem~\cite{Thur} implies that every combinatorial type of
  $3$-polytope has a realization with all edges tangent to the unit
  sphere $S^2$.  Furthermore, the representation is unique up to
  M\"obius transformations; thus, in particular, symmetric
  graphs/lattices have symmetric realizations.  (See Bobenko \&
  Springborn~\cite{BobenkoSpringborn} for a powerful treatment of this
  result, and for references to its involved history.)
\end{compactenum}
This is the situation in dimension~$3$.  The picture in dimension $4$
is not only quite incomplete; it is also clear by now that the results
for the case of $4$-dimensional polytopes will be much more involved.
So, it will be a substantial challenge for 2006 (the centennial of
Steinitz' little 2\hbox{$\,^1\!/_2$}-page paper~\cite{Stei3}) to
characterize the closures of the convex cones of $f$-vectors and
flag-vectors for $4$-polytopes.  This paper sketches some obstacles on
the way as well as some efforts that have been undertaken or that
should and will be made towards this goal. The obstacles are closely
linked to the three results listed above:
\begin{compactenum}[~1.]
\item The geometry of the set of $f$-vectors of $4$-polytopes is rather
  intricate. It does not consist of all the integer points in its
  convex hull, its convex hull is not closed,
  and the cone it spans may be not polyhedral.
\item Moreover, the hierarchy covered by the second Steinitz theorem
  --- connected Eulerian lattices, strongly regular CW spheres, convex
  polytopes, rational convex polytopes --- does not collapse in
  dimension $4$: The set of combinatorial types becomes increasingly
  restricted in this sequence.
\item Furthermore, the non-existence of edge-tangent representations
  for many types of $4$-polytopes is an obstruction to the
  ``E-construction'' (see Section~6) that has recently produced
  sequences of interesting examples.
\end{compactenum}
Nevertheless, there is hope: The
boundary complex of a $4$-dimensional polytope is
$3$-dimensional --- thus we are in essence concerned with
problems of $3$-dimensional combinatorial geometry.
That is, $4$-dimensional polytopes and their faces can
be effectively constructed, handled, and visualized.
The tools that we have available in this context include
dimensional analogy, Schlegel diagrams (see~\cite[Lect.~5]{Z35}),
a connection to tilings that will be
outlined in Section~7, as well as computational tools
(use \url{Polymake}~\cite{GawrilowJoswig}).

\section{{\boldmath{$f$}}-Vectors and flag-vectors} 
\setzero\vskip-5mm \hspace{5mm }

The combinatorial type of a convex $d$-dimensional polytope $P$
(\emph{$d$-polytope}, for short)
is given by its \emph{face lattice} $L(P)$:
This is a finite graded lattice of length~$d+1$
which is \emph{Eulerian}, that is,
every non-trivial interval contains the same number
of elements of odd and of even rank (cf.\ Stanley~\cite{St,Sta5}).
Furthermore, for $d>1$ this lattice is \emph{connected},
that is, the bipartite graph of atoms and coatoms (elements of rank~$1$
vs.\ elements of rank~$d$, corresponding to vertices vs.\ facets) is connected.

The primary numerical data of a polytope, or much more generally
of a graded lattice, are the numbers $f_i=f_i(P)$ of $i$-dimensional
faces ($i$-faces, resp.\ lattice elements of rank~$d+1$).
More generally, one considers the $2^d$
flag numbers $f_S=f_S(P)$ (for $S\subseteq\{0,1,\dots,d-1\}$)
that count the chains of faces with one $i$-face
 for each $i\in S$. These are collected to yield the
\emph{$f$-vector} $f(P):=(f_0,f_1,\dots,f_{d-1})\in\Z^d$,
and the \emph{flag-vector} $\flag(P):=(f_S: S\subseteq\{0,1,\dots,d-1\})$
of the polytope or lattice.
In terms of flag numbers, the bipartite graph of atoms and coatoms
has $f_0+f_{d-1}$ vertices and $f_{0,d-1}$ edges.

For general polytopes, the components of the $f$-vector satisfy
only one non-trivial linear equation, the Euler-Poincar\'e
relation $f_0-f_1\pm\dots+(-1)^{d-1}f_{d-1}=1+(-1)^{d-1}$. The
flag-vector (with $2^d$ components, including the $f$-vector) is
highly redundant, due to the linear ``generalized Dehn-Sommerville
relations'' (Bayer \& Billera~\cite{BaBi}) that allow one to
reduce the number of independent components to $F_d-1$, one less
than a Fibonacci number. In particular, for $d=3$ there is no
additional information in the flag-vector, by
$f_{01}=f_{12}=2f_1$, $2f_{02}=f_{012}=4f_1$. For $4$-polytopes,
the full flag-vector is determined by
\[
\fourflag(P):=(f_0,f_1,f_2,f_3;f_{03}).
\]
(We do \emph{not} delete one of the $f_i$ via the Euler-Poincar\'e
relation, in order to explicitly retain the symmetry for dual
polytopes.) As an example, the flag-vector of the $4$-simplex is
given by $\fourflag(\Delta_4)=(5,10,10,5;20)$. The set of all
$f$-resp.\ flag-vectors of $4$-polytopes will be denoted
$\fvectors(\PP)$ resp.\ $\flagvectors(\PP)$.

The known facts about and partial description of the sets $\fvectors(\PP)$
and $\flagvectors(\PP)$ have been reviewed in detail in
Gr\"unbaum~\cite[Sect.~10.4]{Gr1},
Bayer~\cite{Bay}, Bayer \& Lee~\cite[Sect.~3.8]{BaLee}, and
H\"oppner \& Ziegler~\cite{Z59}.
Here we will be concerned only with
the \emph{linear} known conditions that are tight at
$\fourflag(\Delta_4)$, and
concentrate of the case of $f$-vectors rather than flag-vectors.

\section{Geometry/Topology/Combinatorics} 
\setzero\vskip-5mm \hspace{5mm }

It pays off to study the $f$- and flag-vector problems with
respect to the following hierarchy of four models---a
combinatorial, a topological, and two geometric ones (where the
last one includes arithmetic aspects):
\begin{compactdesc}
\item[Eulerian lattices:]
   Let $\EL$ be the class of all connected Eulerian lattices of
   length~$5$, as defined/reviewed above.
   (More restrictively, one could require that all intervals
   of length at least~$3$ must be connected.)
\item[Cellular spheres:] Let $\SP$ be the class strongly regular
  cellulations of the $3$-sphere, that is, of all regular cell
  decompositions of~$S^3$ for which any intersection of two
  cells is a face of both of them (which may be empty). These objects
  appear as ``regular CW $3$-spheres with the intersection property''
  as in Bj\"orner~et al.\ \cite[pp.~203, 223]{Z10};
  following Eppstein, Kuperberg \&
  Ziegler~\cite{Z80} we call them ``strongly regular spheres.'' The
  intersection property is equivalent to the fact that the face poset
  of the cell complex is a lattice.
\item[Convex polytopes:]
$\PP$ denotes the combinatorial types of convex $4$-polytopes.
\item[Rational convex polytopes:]
$\PR$ will denote the combinatorial types of
convex $4$-polytopes that have a realization with rational
(vertex) coordinates.
\end{compactdesc}
We have natural inclusions
\[
\PR\ \ \subset\ \ \PP\ \ \subset\ \ \SP\ \ \subset\ \ \EL.
\]
The first inclusion is strict due to the existence of non-rational
$4$-polytopes (Richter-Gebert~\cite{Rich4}), the second one due to the
known examples of non-realizable triangulated $3$-spheres,
the third since any strongly regular cell decomposition (e.\,g., a
triangulation) of a compact connected $3$-manifold without boundary
has a connected Eulerian face lattice of length~$5$.

For each of the four classes we define its \emph{cone of flag-vectors},
that is, the closure of the cone with apex~$\fourflag(\Delta_4)$ that is
spanned by all vectors of the form $\fourflag(P)-\fourflag(\Delta_4)$:
For each family of combinatorial types
we denote by $\flagvectors(\cdot)$ the set of flag-vectors, and
by $\closedfcone(\cdot)$ resp.\ $\closedflagcone(\cdot)$
the corresponding closures of the cones of $f$-vectors resp.\
flag-vectors.
So we get the inclusions
\[
\fvectors(\PR)\ \subseteq\
\fvectors(\PP)\ \subseteq\
\fvectors(\SP)\ \subseteq\
\fvectors(\EL),
\]
and
\[
\closedfcone(\PR)\ \ \subseteq\ \
\closedfcone(\PP)\ \ \subseteq\ \
\closedfcone(\SP)\ \ \subseteq\ \
\closedfcone(\EL),
\]
and similarly for flag-vectors ---
\emph{but can we separate them? Is any of these inclusions strict?}
At the moment, that does not appear to be clear, not even if we
consider the sets of $f$- or flag-vectors themselves rather than just the
closures of the cones they span!

\section{Fatness and complexity} 
\setzero\vskip-5mm \hspace{5mm }

Instead of linear combinations of face numbers or flag numbers (such
as the toric $h$-vector, the $cd$-index etc.~\cite{BaLee}), in the
following we will rely on quotients of such. Thus we obtain
homogeneous ``density'' parameters that characterize extremal
polytopes. Such a quotient is the \emph{average vertex degree}
$\frac{f_{01}}{f_0}=2\frac{f_1}{f_0}$.  However, we will in addition
normalize the quotients such that numerator and denominator vanish for
the simplex; then every inequality of the form
``our~parameter$\,\ge\,$constant'' translates into a linear inequality
that holds with equality at the simplex.  Thus instead of the average
vertex degree we would use densities like
$\delta_0:=\frac{f_1-10}{f_0-5}$ or
$\delta_0:=\frac{f_1+2f_3-20}{f_0+f_3-10}$.  (They are not defined for
the simplex.)  For both of these densities
equality in the valid inequality $\delta_0\ge2$
characterizes simple $4$-polytopes.

We prefer such density parameters
since they provide measures of complexity that are independent
of the (combinatorial) ``size'' of the polytope.
For example, they are essentially stable under
various operations of ``glueing'' polytopes or ``connected sums''
of polytopes (as in~\cite[p.~274]{Z35}). In terms of the closed cones
of $f$- and flag-vectors, the
density parameters measure ``how close we are to the boundary''
in terms of a ``slope.''
The following two parameters we call \emph{fatness} and \emph{complexity}:
\[
\fatness(P)\ :=\ \frac{f_1+f_2-20}{f_0+f_3-10},
\qquad
\complexity(P)\ :=\ \frac{\ \ f_{03}\ \ \ -20 }{f_0+f_3-10}.
\]
Both parameters are self-dual: Any polytope and its polar dual
have the same fatness and complexity.
The definition of fatness given here differs from that in~\cite{Z80} by
the additional normalization, which makes the inequalities below
come out simpler. Using the generalized Dehn-Sommerville equations,
one derives from $f_{023}\ge 3f_{03}$ resp.\ from
$f_{02}\ge3f_2$ and $f_{13}\ge3f_1$ that
\[
\complexity(P)\ \le\ 2\fatness(P)-2\qquad\textrm{and}\qquad
\fatness(P)\ \le\ 2\complexity(P)-2,
\]
with equality in the first case if and only if all facets of~$P$ are
simple, and with equality in the second case if and only if
$P$ is $2$-simple and $2$-simplicial.
In particular, $\complexity(P),\fatness(P)\ge2$ holds for all
polytopes, and more generally for all Eulerian lattices of length~$5$.
Furthermore we see that the two parameters are asymptotically
equivalent: Any polytope of high fatness has high complexity,
and conversely.

In terms of complexity and fatness, we can e.\,g.\ rewrite the
flag-vector inequality conjectured by Bayer~\cite[p.~145]{Bay}
as $\fatness(P)\ge 2\complexity(P)-5$;
counter-examples appear in Section~6.
Similarly, the conjectured $f$-vector inequality of Bayer~\cite[p.~149]{Bay}
becomes self-dual if we add the dual inequality  ---  then
we obtain the condition $\fatness(P)\le5$, which again
fails on examples given below.

Another important, self-dual inequality reads
$f_{03}-3(f_0+f_3)\ge-10$, that is,
\[
\complexity(P)\ \ge\ 3.
\]
In fact, this is the $4$-dimensional case of the condition
``$g^{\mathrm {tor}}_2(P)\ge0$'' on the toric $h$-vector~\cite{Sta7}
of a polytope; it was proved by Stanley for rational polytopes only,
and verified for all convex polytopes by Kalai~\cite{kalai87:_rigid_i}
with a simpler argument based on a rigidity result of Whiteley.
It is still not established for the case of strongly regular
$3$-spheres, or for Eulerian lattices of length~$5$.


\section{The \boldmath{$f$}-Vector cone for \boldmath{$4$}-polytopes}
\setzero\vskip-5mm \hspace{5mm }

\noindent
{\bf Theorem} (cf.\ Bayer~\cite[Sect.~4]{Bay}, Eppstein et
al.~\cite{Z80}){\bf.}\\
\it
In terms of the homogeneous coordinates
$\varphi_0:=\frac{f_0-5}{f_1+f_2-20}$ and
$\varphi_3:=\frac{f_3-5}{f_1+f_2-20}$,
for $(f_1+f_2-20,f_0-5,f_3-0)$-space,
the five linear inequalities\rm
\begin{compactenum}[\rm (i) ]
\item $ \varphi_0+3\varphi_3\le1$ \quad (with equality for simplicial
polytopes),
\item $3\varphi_0+ \varphi_3\le1$ \quad (with equality for simple
polytopes),
\item $ \varphi_0\ge0$ \quad (close-to-equality if there are
                              much fewer vertices than other
                              faces),
\item $ \varphi_3\ge0$ \quad (close-to-equality if there are
                              much fewer facets   than other
                              faces),
\item $ \varphi_0+\varphi_3\le\frac25$ \quad (with equality if
           $\fatness(P)=\frac1{\varphi_0+\varphi_3}=\frac52$, which forces
  $g^{\mathrm {tor}}_2(P)=0$)
\end{compactenum}
\it define a $3$-dimensional closed polyhedral cone with apex
$\fourflag(\Delta_4)=(0,0,0)$. It is a cone over a pentagon, which
drawn to scale looks as follows:

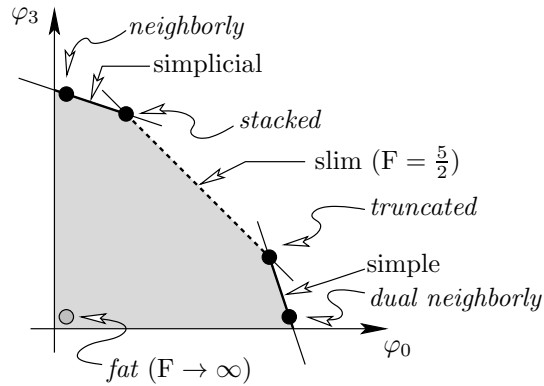
\centerline{\begin{picture}(0,0)%
\includegraphics{fourcone2.pstex}%
\end{picture}%
\setlength{\unitlength}{3947sp}%
\begingroup\makeatletter\ifx\SetFigFont\undefined%
\gdef\SetFigFont#1#2#3#4#5{%
  \reset@font\fontsize{#1}{#2pt}%
  \fontfamily{#3}\fontseries{#4}\fontshape{#5}%
  \selectfont}%
\fi\endgroup%
\begin{picture}(2349,2415)(964,-3106)
\put(1126,-811){\makebox(0,0)[rb]{\smash{\SetFigFont{10}{12.0}{\familydefault}{\mddefault}{\updefault}{\color[rgb]{0,0,0}$\varphi_3$}%
}}}
\put(1676,-906){\makebox(0,0)[lb]{\smash{\SetFigFont{10}{12.0}{\familydefault}{\mddefault}{\updefault}{\color[rgb]{0,0,0}\textit{neighborly}}%
}}}
\put(1836,-1141){\makebox(0,0)[lb]{\smash{\SetFigFont{10}{12.0}{\familydefault}{\mddefault}{\updefault}{\color[rgb]{0,0,0}\textrm{simplicial}}%
}}}
\put(2371,-1496){\makebox(0,0)[lb]{\smash{\SetFigFont{10}{12.0}{\familydefault}{\mddefault}{\updefault}{\color[rgb]{0,0,0}\textit{stacked}}%
}}}
\put(2851,-1751){\makebox(0,0)[lb]{\smash{\SetFigFont{10}{12.0}{\familydefault}{\mddefault}{\updefault}{\color[rgb]{0,0,0}\textrm{slim ($\fatness=\frac52$)}}%
}}}
\put(3191,-2041){\makebox(0,0)[lb]{\smash{\SetFigFont{10}{12.0}{\familydefault}{\mddefault}{\updefault}{\color[rgb]{0,0,0}\textit{truncated}}%
}}}
\put(3176,-2401){\makebox(0,0)[lb]{\smash{\SetFigFont{10}{12.0}{\familydefault}{\mddefault}{\updefault}{\color[rgb]{0,0,0}\textrm{simple}}%
}}}
\put(3186,-2611){\makebox(0,0)[lb]{\smash{\SetFigFont{10}{12.0}{\familydefault}{\mddefault}{\updefault}{\color[rgb]{0,0,0}\textit{dual neighborly}}%
}}}
\put(3281,-2901){\makebox(0,0)[lb]{\smash{\SetFigFont{10}{12.0}{\familydefault}{\mddefault}{\updefault}{\color[rgb]{0,0,0}$\varphi_0$}%
}}}
\put(1546,-3061){\makebox(0,0)[lb]{\smash{\SetFigFont{10}{12.0}{\familydefault}{\mddefault}{\updefault}{\color[rgb]{0,0,0}\textit{fat \rm($\fatness\rightarrow\infty$)}}%
}}}
\end{picture}
}

The five linear inequalities are
valid for $\closedfcone(\PP)$, and they are
tight and facet-defining for $\closedfcone(\SP)$.
\rm

The interesting/crucial parts of this theorem are
on the one hand the existence of arbitrarily fat objects,
which was established for strongly cellular spheres
in~\cite[Sect.~4]{Z80} by the $M_g$-construction (see Section~6),
and the last inequality (v), which follows for polytopes from
$\complexity(P)\ge3$~\cite{Bay}, but whose validity for
strongly regular spheres is an open problem.

Thus \emph{if} $\fatness(P)\ge\frac52$ is valid for
all strongly regular $3$-spheres, then the five inequalities
above give a \emph{complete} linear description of $\closedfcone(\SP)$.
On the other hand, \emph{if} fatness is not bounded for (rational) convex
$4$-polytopes, then the above system is a \emph{complete} description
of~$\closedfcone(\PP)$ resp.~$\closedfcone(\PR)$.
But if one of the two \emph{if}s fails, then the two cones of $f$-vectors
differ substantially!

One can attempt to give a similar description for the $4$-dimensional
cones $\closedflagcone(\PP)$ and $\closedflagcone(\SP)$.
However, in this case the picture (compare Bayer~\cite[Sect.~2]{Bay}
and H\"oppner \& Ziegler~\cite{Z59})
is much less complete, yet.

In both the $f$-vector and in the flag-vector case
the simple polytopes and the simplicial polytopes
appear as extreme cases, and they induce facets that meet
only at the apex (the simplex). The $f$- and flag-vectors
of simple/simplicial $4$-polytopes and $3$-spheres are
well-known --- a complete picture is given by
the $g$-Theorem (McMullen~\cite{McM1}), which for $4$-polytopes
was first established by Barnette~\cite{barnette72:_inequal}
and for $3$-spheres by Walkup~\cite{walkup70}).

\section{Constructions} 
\setzero\vskip-5mm \hspace{5mm }

In order to prove completeness for linear descriptions of
$f$- or flag-vector cones, one needs to have at one's disposal enough
examples or construction techniques for
extremal polytopes that go beyond the usual classes of ``simple and
simplicial'' polytopes (neighborly, stacked, random, etc.).

\para{Cubical polytopes.} (all of whose proper faces are
combinatorial cubes) form a natural class of polytopes.
A very specific construction by Joswig \& Ziegler~\cite{Z62} produced
``neighborly cubical'' polytopes as special projections of
suitably deformed $n$-cubes to~$\R^4$: These are
\emph{rational} cubical $4$-polytopes $C_4^n$
with the graph of the $n$-cube (for $n\ge4$), hence with flag-vectors
\[
\fourflag(C_4^n)\ =\ \big(4,2n,3(n-2),n-2;8(n-2)\big)\cdot 2^{n-2}.
\]
Thus we have rational $4$-polytopes of fatness $\fatness(C_4^n)$
arbitrarily close to~$5$, and complexity $\complexity(C_4^n)$
arbitrarily close to~$8$.  (One may also show that all cubical
polytopes and spheres satisfy $\fatness(P)<5$ and $\complexity<8$.
Indeed, for a polytope of very high fatness and complexity, the facets
need to have a very high number of vertices on average).

\para{The E-construction.}
Eppstein, Kuperberg \& Ziegler~\cite{Z80} presented and analyzed a
particular $4$-dimensional construction that produces
interesting example polytopes:
Let $P\subset\R^4$ be a simple $4$-polytope whose ridges ($2$-faces)
are tangent to~$S^3$; then its polar dual $P^\Delta$ is a
simplicial edge-tangent polytope.
The \emph{E-polytope} of~$P$, obtained as
$E(P):=\conv(P\cup P^\Delta)$, is then $2$-simple and $2$-simplicial.
It has fatness
\[
\fatness(E(P))\ =\ \frac{6f_0(P)-10}{f_0(P)+f_3(P)-5}\ <\ 6.
\]
In~\cite{Z80}, this construction was used to
produce infinite families of $2$-simple $2$-simplicial
polytopes --- apparently the first of their kind.
It was also used to construct $4$-polytopes of fatness
larger than $5.048$ --- currently this is the largest
value that has been observed for convex polytopes.
(All simple and simplicial polytopes have fatness smaller than~$3$.)
We note, however, that the prerequisites for the
E-construction are rather hard to satisfy.
Obvious examples where they can be achieved arise from regular convex
polytopes. On the other hand,
it turned out that, for example,
$P^\Delta$ cannot be a stacked $4$-polytope with more than $6$
vertices! Moreover, for most examples the tangency-condition
seems to force non-rational coordinates for $P$, and hence for $P^\Delta$.
The analysis in~\cite{Z80} depends on a geometric analysis
that puts the Klein model of hyperbolic geometry onto the interior
of the $4$-ball bounded by~$S^3$.

Thus one may ask: Does the E-construction produce non-rational
polytopes? Are there possibly flag-vectors of $2$-simple $2$-simplicial
polytopes that cannot be realized by rational polytopes?
While it seems quite reasonable that the
E-polytope of the regular $120$-cell, with flag-vector
$\fourflag(E(P_{120}))=(720,3600,3600,720;5040)$,
fatness $\fatness(E(P_{120}))>5.02$,
and $720$ biyramids over pentagons as facets, could be non-rational,
current investigations (Paffenholz~\cite{Paffenholz-pc})
suggest that E-polytopes are less rigid than one would
think at first glance, since in some cases the tangency conditions may be
relaxed or dropped.

\para{Fat 3-spheres: The \boldmath{$M_g$}-construction.}
Based on a covering space argument,
Eppstein, Kuperberg \& Ziegler~\cite[Sect.~4]{Z80} constructed
a family of strongly regular CW $3$-spheres whose fatness
is not bounded at all. The construction starts with
a perfect cellulation of the compact connected orientable
$2$-manifold $M_g$ of genus~$g$ with $f$-vector $(1,2g,1)$ and ``fatness''
$\frac{f_1}{f_0+f_2}=g$. Then one shows that there is a
finite covering $\widetilde M_g$ of $M_g$ whose induced cell decomposition
(of the same ``fatness''~$g$) is strongly regular.
Finally, from the standard embedding of $\widetilde M_g\times I$
into~$S^3$, where the interval~$I$ is subdivided very finely
and $\widetilde M_g\times I$ gets the product decomposition,
one obtains a cellulation of $S^3$ whose flag vector
is dominated by the flag-vector of $\widetilde M_g\times I$.
This yields strongly regular cell decompositions of~$S^3$
whose $f$-vector is approximately proportional to
$(1,2g,1)*(1,1)=(1,2g+1,2g+1,1)$.
Thus the resulting spheres have fatness arbitrarily close to $2g+1$.

\para{Many triangulated \boldmath$3$-spheres.}
Applied to the fat $3$-spheres produced by the
$M_g$-construction, the E-construction yields $3$-spheres
with substantially more non-simplicial facets than their number of vertices.
Thus one obtains (Pfeif{}le \cite{Pfeifle-pc}) that on a large
number of vertices there are far more triangulated $3$-spheres
than there are types of simplicial $4$-polytopes, thus resolving
a problem of Kalai~\cite{Ka9}.

\section{Tilings} 
\setzero\vskip-5mm \hspace{5mm }

There are close connections between $d$-polytopes
(``polyhedral tilings of~$S^{d-1}$'')
and normal polyhedral tilings of~$\R^{d-1}$. In particular,
from $4$-polytopes one may construct $3$-dimensional tilings,
for example by starting with a regular tiling of $\R^3$
by congruent tetrahedra and then
replacing the tiles by Schlegel diagrams based on a simplex facet.
(The converse direction, from tilings of~$\R^3$ to $4$-polytopes,
is non-trivial: It hinges on
non-trivial liftability restrictions; see Rybnikov~\cite{rybnikov99:_stres}.)

\emph{Normal tilings} are face-to-face tilings of $\R^{d-1}$ by convex
polytopes for which the inradii and circumradii of tiles are bounded
from below resp.\ from above --- see
Gr\"unbaum \& Shephard~\cite[Sect.~3.2]{gruenbaum87:_tilin_patter}.
Of course, all components of an
$f$-vector for tilings would be infinite, but one can try to define
ratios, e.\,g.\ try to find the ``average'' number of vertices per tile.

\para{The Euler formula for tilings.}
Thus, for $\rho>0$, let $f_i(\rho)$ be the number of all faces
of the tiling that intersect the interior of the ball $B^d(\rho)$
of radius~$\rho$ around the origin.
This yields a regular decomposition of an open $d$-ball into convex cells;
via one-point-compactification (e.\,g.~generated by
stereographic projection) by one additional vertex
we obtain a regular cell-decomposition of a $d$-sphere;
thus we obtain~\cite{Wassmer-pc} that for all $\rho>0$
\[
f_0(\rho) - f_1(\rho) +\dots+(-1)^d f_{d-1}(\rho)\ \ =\ \ (-1)^{d-1}.
\]
In particular, this implies that if the limits
$\phi_i:=\lim_{\rho\rightarrow\infty}\frac{f_i(\rho)}{\sum_jf_j(\rho)}$
exist, then they satisfy $0\le\phi_i\le\frac12$.
Furthermore, the existence of the limits~$\phi_i$ is automatic
for tilings that are invariant under a full-dimensional lattice of
trans\-lations, such as the ``tilings by Schlegel diagrams''
suggested above.
In this case the limits $\phi_i$ are strictly positive,
and they satisfy the \emph{Euler formula for tilings},
\[
\phi_0-\phi_1+\phi_2\pm\dots+(-1)^{d-1}\phi_{d-1}\ \ =\ \ 0.
\]
For such tilings, we can also define flag-numbers $\phi_S(\rho)$. Then
the limit $\rho\rightarrow\infty$ of any quotient such as
$\fatness(\rho):=\frac{f_1(\rho)+f_2(\rho)}{f_0(\rho)+f_3(\rho)}$
exists, it is positive and finite, and it coincides with
\[
\fatness({\cal T})\ \ :=\ \
\frac{\phi_1+\phi_2}{\phi_0+\phi_3}.
\]
One can construct ``tilings by Schlegel diagrams''
with a high (average) number of vertices per tile,
or with a high number of tiles at each vertex.
But are both achievable simultaneously? Equivalently,
is there is a uniform upper
bound on the fatness $\fatness({\cal T})$ of normal tilings
of~$\R^3$? This is not known, yet, but if fatness is bounded for normal
$3$-tilings, then it is bounded for $4$-polytopes as well.

\para{Fat tilings.}
Remarkably, there are tilings that have considerably larger fatness
than the fattest polytopes we know.  In particular,
a modified E-construction applied to suitable Schlegel
$3$-diagrams~\cite[Lect.~5]{Z35} of $C_n\times C_n$ and embedding
into a cubic tiling one obtains normal lattice-transitive tilings of
fatness arbitrarily close to~$6$.

\para{Acknowledgements.}
This report has benefitted from joint work and many fruitful
discussions with David Eppstein, Michael Joswig,
Greg Kuperberg, Julian Pfeif{}le, Andreas Paffenholz, and Arnold Wa{\ss}mer.
(And thanks to Torsten Heldmann!)


\baselineskip 4.1mm

\begin{thebibliography}{10}

\bibitem{barnette72:_inequal}
{D.~W. Barnette}, {\em Inequalities for f-vectors of $4$-polytopes},
  Israel J. Math. {\bf 11} (1972), 284--291.

\bibitem{Bay}
{M.~M. Bayer}, {\em The extended f-vectors of $4$-polytopes}, J.
  Combinat. Theory, Ser.~A {\bf 44} (1987), 141--151.

\bibitem{BaBi}
{M.~M. Bayer and L.~J. Billera}, {\em Generalized {D}ehn-{S}ommerville
  relations for polytopes, spheres and {E}ulerian partially ordered sets},
  Inventiones Math. {\bf 79} (1985), 143--157.

\bibitem{BaLee}
{M.~M. Bayer and C.~W. Lee}, {\em Combinatorial aspects of convex
  polytopes}, in: Handbook of Convex Geometry (P.~Gruber, J.~Wills, eds.),
  North-Holland, Amsterdam 1993, 485--534.

\bibitem{Z10}
    {A.~Bj\"orner, M.~{Las Vergnas}, B.~Sturmfels, N.~White,
    G.~M.  Ziegler}, {\em Oriented Matroids}, {\sl Encyclopedia of
    Math.} {\bf 46}, Cambridge Univ.\ Press, 2nd~ed., 1999.

\bibitem{BobenkoSpringborn}
{A.~I. Bobenko and B.~A. Springborn}, {\em Variational principles for
  circle patterns, and {K}oebe's theorem}, Report No. 545 (Sfb 288
  preprint series), TU~Berlin 2002;
\newblock arXiv: math.GT/0203250, 38.

\bibitem{Z80}
{D.~Eppstein, G.~Kuperberg, and G.~M. Ziegler}, {\em Fat $4$-polytopes and
  fatter $3$-spheres}.
\newblock Preprint, March 2002, 12; arXiv: math.CO/0204007;
\newblock to appear in: W. Kuperberg Festschrift (A. Bezdek, ed.),
                  Marcel-Dekker 2002.

\bibitem{GawrilowJoswig}
{E.~Gawrilow and M.~Joswig}, {\em Polymake: {A} software package for
  analyzing convex polytopes}.
\newblock http://www.math.tu-berlin.de/diskregeom/polymake/

\bibitem{Gr1}
{B.~Gr\"unbaum}, {\em Convex {P}olytopes}, Interscience, London, 1967.
\newblock Revised edition (V. Kaibel, V. Klee, G. M. Ziegler, eds.),
  Springer-Verlag 2002, in preparation.

\bibitem{gruenbaum87:_tilin_patter}
{B.~Gr\"unbaum and G.~C. Shephard}, {\em Tilings and Patterns}, Freeman,
  New York, 1987.

\bibitem{Z59}
{A.~H\"oppner and G.~M. Ziegler}, {\em A census of flag-vectors of
  $4$-polytopes}, in: Polytopes --- Combinatorics and Computation
  (G.~Kalai, G.~Ziegler, eds.), {\sl DMV Seminars} {\bf 29},
  Birkh\"auser-Verlag, Basel 2000, 105--110.

\bibitem{Z62}
{M.~Joswig and G.~M. Ziegler}, {\em Neighborly cubical polytopes}, Discrete
  Comput. Geometry {\bf 24} (2000), 325--344.

\bibitem{kalai87:_rigid_i}
{G.~Kalai}, {\em Rigidity and the lower bound theorem, {I}}, Inventiones
  Math. {\bf 88} (1987), 125--151.

\bibitem{Ka9}
\leavevmode\vrule height 2pt depth -1.6pt width 23pt,
  {\em Many triangulated spheres}, Discrete Comput. Geometry {\bf 3}
  (1988), 1--14.

\bibitem{McM1}
{P.~McMullen}, {\em The numbers of faces of simplicial polytopes}, Israel
  J. Math. {\bf 9} (1971), 559--570.

\bibitem{Paffenholz-pc}
{A.~Paffenholz}, {\em Work in progress}.
\newblock TU Berlin, 2002.

\bibitem{Pfeifle-pc}
{J.~Pfeifle}, {\em Work in progress}.
\newblock TU Berlin, 2002.

\bibitem{Rich4}
{J.~Richter-Gebert}, {\em Realization Spaces of Polytopes},
  {\sl Lecture Notes in Mathematics} {\bf 1643}, Springer-Verlag,
  Berlin Heidelberg, 1996.

\bibitem{rybnikov99:_stres}
{K.~Rybnikov}, {\em Stresses and liftings of cell complexes}, Discrete
  Comput. Geometry {\bf 21} (1999), 481--517.

\bibitem{Schla}
{L.~Schl\"afli}, {\em Theorie der vielfachen Kontinuit\"at}, Denkschriften
  der Schwei\-ze\-ri\-schen naturforschenden Gesellschaft, Vol.~38, 1--237,
  Z\"urcher und Furrer, Z\"urich 1901 (written 1850-1852).

\bibitem{Sta7}
{R.~P. Stanley}, {\em Generalized $h$-vectors, intersection cohomology of
  toric varieties, and related results}, in: Commutative Algebra and
  Combinatorics (M. Nagata, H.~Matsumura, eds.), {\sl Advanced Studies
  in Pure Mathematics} {\bf 11}, Kinokuniya, Tokyo 1987, 187--213.

\bibitem{St}
\leavevmode\vrule height 2pt depth -1.6pt width 23pt, {\em A survey of
  {E}ulerian posets}, in:\ ``Polytopes:
  Abstract, convex and computational'' (Toronto 1993),
  Proc.\ NATO Advanced Study Institute (T.~Bisztriczky et al., eds.),
  Kluwer, Dordrecht 1994, 301--333.

\bibitem{Sta5}
\leavevmode\vrule height 2pt depth -1.6pt width 23pt, {\em Enumerative
  Combinatorics, Volume I}, Second edition, Cambridge Studies in Advanced
  Mathematics {\bf 49}, Cambridge University Press 1997.

\bibitem{Stei3}
{E.~Steinitz}, {\em {\"U}ber die {E}ulerschen {P}olyederrelationen}, Archiv
  f\"ur Mathematik und Physik {\bf 11} (1906), 86--88.

\bibitem{Stei1}
\leavevmode\vrule height 2pt depth -1.6pt width 23pt, {\em {P}olyeder und
  {R}aumeinteilungen}, in Encyklop\"adie der mathematischen
  Wis\-sen\-schaf\-ten,
  Geometrie, III.1.2., Heft~9, Kapitel~\mbox{III\,A\,B\,12},
  (W.~F. Meyer and H.~Mohrmann, eds.),
  B.~G.~Teubner, Leipzig, 1922, 1--139.

\bibitem{Thur}
{W.~P. Thurston}, {\em Geometry and topology of $3$-manifolds}.
\newblock Lecture Notes, Princeton University, Princeton 1977--1978.

\bibitem{walkup70}
{D.~W. Walkup}, {\em The lower bound conjecture for $3$- and
  $4$-manifolds}, Acta Math. {\bf 125} (1970), 75--107.

\bibitem{Wassmer-pc}
{A.~Wa{\ss}mer}, {\em Work in progress}.
\newblock TU Berlin, 2002.

\bibitem{Z35}
{G.~M. Ziegler}, {\em Lectures on {P}olytopes}, {\sl Graduate Texts
  in Mathematics} {\bf 152}, Springer-Verlag, New York, 1995.
\newblock Revised edition, 1998; ``Updates, corrections, and more'' at
  www.math.tu-berlin.de/~ziegler.

\label{lastpage}

\end{thebibliography}
\end{document}